\documentclass{amsart}
\usepackage{amssymb,euscript,amsmath, mathrsfs}
\usepackage{tikz}
\usepackage{bm}
\usepackage{rotating}
\usepackage{graphicx}
\newcounter{ENUM}

\newcommand{\margh}[1]{}

\input xy
\xyoption{all}
\CompileMatrices

\def\bgamma{{\bm{\gamma}}}
\def\bbeta{{\bm{\beta}}}

\def\ZZ{{\mathbb Z}}

\def\AA{{\mathbb A}}

\def\PP{{\mathbb P}}

\def\cP{{\mathcal P}}

\def\fp{{\mathfrak p}}

\def\vp{\varphi}

\def\codim{\operatorname{codim}}

\newtheorem{thm}{Theorem}[section]
\newtheorem{prop}[thm]{Proposition}

\newtheorem{cor}[thm]{Corollary}

\theoremstyle{definition}
\newtheorem{defn}[thm]{Definition}

\newtheorem{ex}[thm]{Example}

\theoremstyle{remark}
\newtheorem{notn}[thm]{Notation}
\newtheorem{rem}[thm]{Remark}

\numberwithin{equation}{section}
\begin{document}
\title{Multigraded Cayley-Chow forms}
\author{Brian Osserman}
\author{Matthew Trager}
\begin{abstract} We introduce a theory of multigraded 
Cayley-Chow forms associated to subvarieties of products of projective
spaces. Two new phenomena arise: first, the construction turns out to require 
certain inequalities on the dimensions of projections; and second,
in positive characteristic the multigraded Cayley-Chow forms can have 
higher multiplicities. The theory also provides a natural framework for
understanding multifocal tensors in computer vision. 
\end{abstract}

\thanks{The first author is partially supported by a grant from the Simons
Foundation \#279151. The second author was supported in part by the ERC 
advanced grant VideoWorld, the Institut Universitaire de France, the 
Inria-CMU associated team GAYA, and the ANR grant RECAP.
This project arose out of discussions at the 
American Institute of Mathematics workshop on Algebraic Vision.}

\maketitle

\margh{Nice condition for when there exists beta which works? guess 
hypersurface ineqs can always be satisfied; oops, not w/given ineqs on
$\beta_i$; cf circuit case. other equivalent to multideg
support full-dimensional? not true for general polytopes, but these are
special}
\margh{something about grassmann tensors?}
\margh{Nonlinear cameras?}

\section{Introduction}

Let $X \subseteq \PP^n$ be a projective variety of dimension $r$ and
degree $d$. The set of all linear spaces of dimension $n-r-1$ meeting
$X$ is a hypersurface $Z_X$ in the Grassmannian ${\rm Gr}(n-r-1,n)$. Any such 
hypersurface can be written as the zero set inside the Grassmannian of
a polynomial $H_X$ in the Pl\"ucker coordinates, which turns out to also be of
degree $d$. The polynomial $H_X$ is known as the ``Chow form'' or
``Cayley form'' of $X$, and we will refer to it as the \textbf{Cayley-Chow 
form}. From the Cayley-Chow form $H_X$ we immediately recover the 
hypersurface $Z_X$, and one can then recover $X$ as the set of points
$P \in \PP^n$ such that every $(n-r-1)$-dimensional linear space containing
$P$ corresponds to an element of $Z_X$. Thus, the Cayley-Chow form can
be used to encode subvarieties of projective space, and this classical
construction has played an important role in moduli space theory,
especially in the guise of Chow varieties, but also for instance in 
Grothendieck's original construction of Quot and Hilbert schemes \cite{gr1}.
See \S 2 of Chapter 3 of \cite{g-k-z1} for a presentation of 
this material.

The purpose of this paper is to generalize this classical theory to the
case of subvarieties of products of projective spaces. We find that
the generalization displays some interesting properties, particularly that 
certain
dimensional inequalities have to be satisfied in order for it to work.
We give necessary and sufficient conditions for the generalized theory
to go through. Moreover, we explain that in cases where the necessary
inequalities are not satisfied, the failure of the theory can in fact
shed light on previously observed phenomena in computer vision. In
addition, positive-characteristic phenomena arise in our
more general setting, causing the Cayley-Chow form to sometimes have
higher multiplicities. 

Before we state our main results, we need to recall the notion of
\textbf{multidegree} of a subvariety 
$X \subseteq \PP^{n_1} \times \cdots \times \PP^{n_k}$. 
If $X$ has codimension $c$, we can represent its multidegree
as a homogeneous polynomial
$$\sum_{\bgamma=(\gamma_1,\dots,\gamma_k)} 
a_{\bgamma} t_1^{\gamma_1} \cdots t_k^{\gamma_k}$$
of degree $c$, where $a_{\bgamma}$ is determined as the
number of points of intersection (counting multiplicity) of $X$ with
$L_1 \times \cdots \times L_k$, where each $L_i$ is a general linear
space of dimension $\gamma_i$. The multidegree is also equivalent to the 
Chow class of $X$, although we will not need this.

We can summarize our main results in the following theorem:

\begin{thm}\label{thm:main} Let 
$X \subseteq \PP^{n_1} \times \ldots \times \PP^{n_k}$ be a 
  projective variety of dimension $r$, and suppose that $X$ is not of
the form $X' \times \prod_{i \not\in I} \PP^{n_i}$ for any 
$I \subsetneq \{1,\dots,k\}$ and $X' \subseteq \prod_{i \in I} \PP^{n_i}$.
Given also 
  $\bbeta=(\beta_1,\dots,\beta_k)$ with $0 \leq \beta_i \leq n_i$ for
  $i=1,\dots,k$ and $\sum_{i=1}^k \beta_i = r+1$, write 
  $\alpha_i=n_i-\beta_i$ for each $i$.

  Consider the closed subset
$$Z_{X,\bbeta} 
=\{(L_1,\dots,L_k):X \cap (L_1 \times \cdots \times L_k)\neq \emptyset \} 
\subseteq 
{\rm Gr}(\alpha_1, n_1) \times \cdots \times {\rm Gr}(\alpha_k, n_k).$$ 
Then $Z_{X,\bbeta}$ is a hypersurface determining $X$ 
if and only if for every nonempty 
  $I \subsetneq \{1,\dots, k\}$ we have 
$$\dim p_I(X) \geq \sum_{i \in I} \beta_i,$$
where $p_I(X)$ denotes the 
projection of $X$ onto $\prod_{i \in I} \PP^{n_i}$.

Assuming the above inequalities are satisfied, we have that 
$Z_{X,\bbeta}$ is
the zero set of a single multihomogeneous polynomial $F_X$ under the
product of the Pl\"ucker embeddings of the ${\rm Gr}(\alpha_i,n_i)$.
There is a multiplicity $\epsilon_{X,\bbeta}\geq 1$ such that if we write
$$H_{X,\bbeta}:= F_X^{\epsilon_{X,\bbeta}},$$
and $X$ has multidegree
$$\sum_{\bgamma} 
a_{\bgamma} t_1^{\gamma_1} \cdots t_k^{\gamma_k},$$
then the multidegree of $H_{X,\bbeta}$ (as a multihomogeneous polynomial) 
is given by
$$(a_{\alpha_1+1,\alpha_2,\dots,\alpha_k}, \dots, 
a_{\alpha_1, \dots, \alpha_{k-1},\alpha_k+1}).$$ 

Finally, in characteristic $0$, we always have
$\epsilon_{X,\bbeta}=1$.
\end{thm}

In fact, we first show in Proposition \ref{prop:hypersurface}
that $Z_{X,\bbeta}$ is a hypersurface if and only if the slightly
weaker inequalities \eqref{eq:hypersurface-ineq} are satisfied, and 
then in Proposition \ref{prop:degenerate} that when $Z_{X,\bbeta}$ 
is a hypersurface, it determines $X$ uniquely if and only if the above 
inequalities are satisfied. The condition that $X$ not be a product with
any of the projective spaces is just to simplify the statement; see
Example \ref{ex:decomposable} below. The multiplicity $\epsilon_{X,\bbeta}$
is defined naturally in Definition \ref{def:multichow} as the degree of 
the map from an incidence 
correspondence to $Z_{X,\bbeta}$, and in positive characteristic, 
it may be strictly greater than $1$; see Example \ref{ex:frobenius} below.
The remaining statements of the theorem are proved in Corollary
\ref{cor:multiplucker}, Proposition \ref{prop:degree-1}, and 
Corollary \ref{cor:multidegree}. An interesting aspect of 
Proposition \ref{prop:degenerate} in comparison to the classical case is 
that if we define the set $S_Z$ to consist of points 
$P \in \PP^{n_1} \times \cdots \times \PP^{n_k}$ with 
the property that every $L_1 \times \cdots \times L_k$ containing $P$ 
must lie in $Z_{X,\bbeta}$, then $X \subseteq S_Z$ always, but we
will have in general that $S_Z$ contains additional components.

In addition, while we have stated the fundamental inequalities in
terms of the dimensions of projections of $X$, according to work of
Castillo, Li and Zhang \cite{c-l-z1} this may be reinterpreted with an
equivalent formulation in terms of the multidegree of $X$. We carry this
out in Corollary \ref{cor:multideg-ineqs} below. Thus, while we do require
non-trivial conditions on $X$ in order for our theory to apply, when it
applies it does so uniformly to all $X$ of given multidegree.

It is natural to wonder how our construction compares to applying a
Segre embedding together with the classical construction; we discuss
this briefly in Remark \ref{rem:segre} below.

We next make the following observation.

\begin{rem}\label{rem:too-many} The inequalities given in Theorem 
\ref{thm:main} can only
be satisfied if $k \leq r+1$, since otherwise the set of $i$ such that
$\beta_i \neq 0$ is necessarily proper in $\{1,\dots,k\}$, and will
violate the necessary inequality. 
\end{rem}

We can now explain the relationship to computer vision, and specifically to
the reconstruction of a configuration of unknown cameras. 
A basic model for a (positioned) camera is as a linear projection from
the three-dimensional world to the two-dimensional film/sensor plane,
which we consider as a $\PP^3$ and a $\PP^2$, respectively. Thus,
a $k$-tuple of cameras corresponds to an $k$-tuple of
linear projections, which together induce a rational map
$$\PP^3 \dashrightarrow (\PP^2)^k.$$
The closure of the image of this map is then a three-dimensional
subvariety of $(\PP^2)^k$, which is called the ``multiview variety.''
This can be thought of as describing which $k$-tuples of points in $\PP^2$
could come from a single point in $\PP^3$. Knowing the multiview variety is
equivalent to knowing the camera configuration, at least up to change of
`world coordinates' in $\PP^3$. It is well known in computer vision that
for $k =2,3,4$, there exists a $k$-tensor which determines the camera
configuration, called the ``multifocal tensor.'' It is equally well known
that this construction does not extend to $k > 4$. 

The multifocal tensor is described in terms of incidences with $k$-tuples 
of linear spaces, and in fact this
inspired our construction. On the other hand, we can reinterpret the theory
of multifocal tensors in terms of Theorem \ref{thm:main} as follows. 
First, Aholt, Sturmfels and Thomas showed in Corollary 3.5 of \cite{a-s-t1}
that all the coefficients of the multidegree of a multiview variety are
equal to $1$, so we see that when our Cayley-Chow form construction applies,
the result is a multilinear polynomial in $k$ variable sets, which is to say,
a $k$-tensor. It is routine to check that for $k \leq 4$, our construction
does apply for suitable choice of $\bbeta$, and the Cayley-Chow form 
coming from the multivew variety is precisely 
the multifocal tensor. Conversely, if $k\geq 5$, then
Remark \ref{rem:too-many} implies that no analogous construction exists
for any choice of $\bbeta$. See Examples \ref{ex:multiview} and
\ref{ex:trifocal} for further details. Beyond giving a new point of view
on these known constructions, we also hope that Theorem \ref{thm:main}
will provide new applications in computer vision, in the context of
generalized cameras. Recent work of Ponce, Sturmfels and the second author
\cite{p-s-t1}, and of Escobar and Knutson \cite{e-k1}
develops a theory of configurations of such cameras,
including multidegree-type formulas, and we expect that Theorem 
\ref{thm:main} will provide a generalization of multifocal tensors to this
setting, where the tensors will be replaced with higher-degree forms.

\begin{rem}
A different connection is to the notion of circuit polynomials in matroid
theory, which we now describe. First observe that when all $n_i$ are equal 
to $1$, we must have $\beta_i=1$
for some subset $S$ of $I$ of size $r+1$, and $\beta_i=0$ for $i \not \in S$.
The inequalities of Theorem \ref{thm:main} are never satisfied except in the 
trivial case that $r=k-1$, and $\beta_i=1$ for all $i$, so that 
$H_{X,\bbeta}$ will simply recover the defining polynomial of $X$. 
However, the weaker inequalities \eqref{eq:hypersurface-ineq} will 
be satisfied more generally: specifically, whenever we have 
$\dim p_I(X) = r$ for the $I$ as above.
% direct translation is $$\dim p_{I'}(X) \geq |I \cap I'|-1$$ for all $I'$.
% enough to consider $I' \subseteq I$, then since dim can only drop by 1
% each time drop a coordinate, reduce to case $I'=I$. 
Thus, in this case we can still
define our multigraded Cayley-Chow form, although it will not suffice to
recover $X$.
% in fact, recovers precisely projection of $X$ to coordinates in $S$
This special case of our construction turns out to be connected to algebraic 
matroids. 

Specifically, one approach to algebraic matroids is as follows: in order 
to construct a matroid on $\{1,\dots,k\}$ of rank $r$, choose a variety 
$X \subseteq \AA^k=(\AA^1)^k$ of dimension $r$, with prime ideal 
$\fp \subseteq K[x_1,\dots,x_k]$. Define independent sets by algebraic 
independence of the images of the $x_i$ modulo $\fp$. In this context, 
a circuit $C \subseteq \{1,\dots,k\}$ will have the property that the
closure of the projection of $X$ to $\prod_{i \in C} \AA^1$ has codimension 
$1$, and
hence is cut out by a single polynomial $H$ in $r+1$ variables. This 
polynomial is called the ``circuit polynomial,'' and precisely cuts out 
the closure of the locus of $(P_i)_{i \in C} \in (\AA^1)^{r+1}$ such that 
there exists $(P_i)_{i \not \in C}$ with $(P_1,\dots,P_k) \in X$. See \S 5 of
Kir\'aly-Rosen-Theren \cite{k-r-t1}. Replacing $X$ by
its closure in $(\PP^1)^k$ and $H$ by its multihomogenization, we recover
the Cayley-Chow form construction with $\beta_i=1$ for $i \in C$ and
$\beta_i=0$ otherwise (at least, up to omission of 
$\epsilon_{X,\bbeta}$). 
% could have epsilon>1 by eg taking $k=3$ and $(t,t^p,t^{p^2})$ in char p.
% not sure about char 0. but this construction always gives multiplicity 1,
% even in ideal version described in paper.
% if take projection of closure, get closed subset in (P^1)^r, whose
% restriction to (A^1)^r is original closure of projection -- one containment
% obvious, other follows from irreducibility, fact that original contains open
% dense subset.
% this situation (more generally if every beta_i is either 0 or n_i), 
% cayley-chow form just cuts out projection, so get desired description
\margh{degree behavior?}
\end{rem}

Finally, we mention that the inequalities arising both in Theorem
\ref{thm:main} and in the condition for $Z_{X,\bbeta}$ to be a 
hypersurface
(see Proposition \ref{prop:hypersurface} below) are closely related
to concepts arising in polymatroid theory; see Remark \ref{rem:polymatroid}
below.

\subsection*{Acknowledgements} We would like to thank Bernd Sturmfels
for bringing to our attention various connections to the literature, 
particularly the notion of matroid circuit polynomials.

\subsection*{Conventions} We work throughout over an algebraically closed
field $K$. A variety is always assumed irreduible.

Given $\bbeta=(\beta_1,\dots,\beta_k)\in \ZZ^k$, we will write 
$|\bbeta|:=\sum_{i=1}^k \beta_i$, and for $I \subseteq \{1,\dots,k\}$, 
we will write $\bbeta_I:=(\beta_i)_{i \in I}$, and 
$|\bbeta_I|:=\sum_{i \in I} \beta_i$. We also write 
$I^c:=\{1,\dots,k\} \smallsetminus I$.

\section{Multidegrees and dimensions of projections}

We begin by collecting some background results on the relationship between 
multidegree and dimensions of projections. If we have a subvariety $X$
of $\prod_{i=1}^k \PP^{n_i}$ with multidegree
$$\sum_{\bgamma} 
a_{\bgamma} t_1^{\gamma_1} \cdots t_k^{\gamma_k},$$
we say the \textbf{support} of the multidegree is the set of 
$\bgamma$ for which $a_{\bgamma} \neq 0$.
Note that by definition, this is contained in the subset of 
$\bgamma$ with $\gamma_i \geq 0$ for all $i$, and
$\sum_i \gamma_i=\codim X$.

The main theorem of \cite{c-l-z1} asserts:

\begin{thm}[Castillo-Li-Zhang]\label{thm:support} If 
$$X \subseteq \PP^{n_1} \times \dots \times \PP^{n_k}$$
is irreducible, the support of its multidegree is 
$$\left\{\bgamma: \sum_i (n_i-\gamma_i)= \dim X, \text{ and }
\sum_{i \in I} (n_i-\gamma_i) \leq \dim p_I(X) \quad \forall I \subsetneq
\{1,\dots,k\}\right\},$$
where $p_I$ denotes projection onto the product of the subset of the
$\PP^n$s indexed by $I$.

Moreover, the function $\delta(I)= \dim p_I X$ satisfies the
following conditions: 
\begin{itemize}
\item $\delta(\emptyset)=0$; 
\item for $I \subseteq J$ we have $\delta(I) \leq \delta(J)$;
\item and for any $I,J$ we have $\delta(I \cap J)+\delta(I \cup J) \leq 
\delta(I) +\delta(J)$. 
\end{itemize}
% NB: statement fails if X allowed reducible; may have non-convex support
\end{thm}

In fact, they treat the case that all $n_i$ are equal, but one reduces
immediately to this case by linearly embedding each $\PP^{n_i}$ into
a larger projective space of fixed dimension.
% doesn't affect dims of projections
% if start w/n_i's, re-embed into fixed n, codimension goes from 
% c to c+\sum_i (n-n_i), if have \sum \alpha_i = c, can replace
% with \alpha_i+(n-n_i) for new monomial. previously had linear subspaces
% of dim alpha_i, now will be alpha_i'=alpha_i+(n-n_i). Then
% alpha_i'-n=alpha_i-n_i, and intersecting
% with general linear space of this dimension same as previously
% intersecting with general linear space of dimension alpha_i. if
% choose alpha_i<n-n_i, will not meet, so outside support.
Note that the first part of the theorem is equivalent to saying that
a general choice of $L_i$ of dimension $\gamma_i$ will yield
$X \cap (L_1 \times \cdots \times L_k) \neq \emptyset$ if and only if
$\sum_{i \in I} (n_i-\gamma_i) \leq \dim p_I(X)$ for all $I \subsetneq
\{1,\dots,k\}$. The second part of the theorem connects multidegrees to
polymatroid theory, and says in particular that the function $\delta$
is ``submodular.''

The first part of the theorem implies that the dimensions of 
the projections determine the support of the
multidegree. We see using some standard facts in polymatroid theory
that the converse also holds.

\begin{cor}\label{cor:multdeg-projs} Given 
$$X \subseteq \PP^{n_1} \times \cdots \times \PP^{n_k}$$
irreducible, the data of the support of the multidegree of $X$ is 
equivalent to the data of the dimensions of $p_I(X)$ for all
$I \subseteq \{1,\dots,k\}$.
\end{cor}

\begin{proof} Let $\cP$ be the polytope cut out by the inequalities of
Theorem \ref{thm:support}, and $\bar{\cP}$ the face of $\cP$ cut out by
the hyperplane $\sum_i (n_i-\gamma_i)= \dim X$. Then $\cP$ is known 
as a ``polymatroid'', and $\bar{\cP}$ is the set of ``bases;'' see
\S 1 of \cite{h-h1}. Moreover, from Proposition 1.3 of \cite{h-h1} 
we see that the vertices of $\bar{\cP}$ are integral, and for every
$I$ there is some vertex lying in the corresponding bounding hyperplane
(in their notation, we take any $\pi$ such that $I=\{i_1,\dots,i_{|I|}\}$).
It follows that we can recover the $\dim p_I(X)$ from the integral points
of $\bar{\cP}$. Since Theorem \ref{thm:support} says that the support of 
the multidegree is equal to the set of lattice points in $\bar{\cP}$, 
we conclude that it determines the $\dim p_I(X)$, as desired.
\end{proof}

The following standard fact from polymatroid theory will also be helpful.
Since the proof is quite short, we include it.

\begin{prop}\label{prop:minimal-subset} Suppose we have a function
$\delta$ from subsets of $\{1,\dots,k\}$ to $\ZZ_{\geq 0}$ satisfying
the conditions in Theorem \ref{thm:support}, and write 
$r=\delta(\{1,\dots,k\})$. Suppose also that we are given 
$\bbeta=(\beta_1,\dots,\beta_k) \in (\ZZ_{\geq 0})^k$ with $|\bbeta|=r+1$,
and satisfying that for all $I \subseteq \{1,\dots,k\}$, we have
\begin{equation}\label{eq:1-deficient}
|\bbeta_I| \leq \delta(I)+1.
\end{equation}
Then there exists a nonempty
$J \subseteq \{1,\dots,k\}$ such that for all $I \subseteq \{1,\dots,k\}$,
we have
\begin{equation}\label{eq:hyper-sharp}
|\bbeta_I| = \delta(I)+1
\end{equation} 
if and only if $I \supseteq J$.
\end{prop}

\begin{proof} First note that if such a $J$ exists, it is necessarily 
nonempty, since we have assumed $\delta(\emptyset)=0$.
Now, if we have $I_1$ and $I_2$ satisfying \eqref{eq:hyper-sharp}, then we
see that
\begin{align*}
\delta(I_1 \cap I_2) 
& \leq \delta (I_1)+\delta (I_2)-\delta (I_1 \cup I_2) \\
& \leq |\beta_{I_1}|-1
+ |\bbeta_{I_2}|-1
- |\bbeta_{I_1 \cup I_2}|-1 \\
& = |\bbeta_{I_1 \cap I_2}|-1,
\end{align*}
so $I_1 \cap I_2$ also satisfies \eqref{eq:hyper-sharp}. The result
follows.
\end{proof}

\begin{rem}\label{rem:polymatroid}
In polymatroid theory, there is a notion of ``$1$-deficient'' vectors,
and among vectors $\bbeta=(\beta_1,\dots,\beta_k)$ satisfying the conditions
that $|\bbeta|=r+1$, those satisfying \eqref{eq:1-deficient} are 
precisely the $1$-deficient vectors. This is exactly the condition that 
arises for us in order for $Z_{X,\bbeta}$ to be a hypersurface -- see 
Proposition \ref{prop:hypersurface}
below. 
% \beta_i \leq n_i+1 follows, but not beta_i \leq n_i. but still correct as
% stated.
Now, given such a $\bbeta$, note that the $J$ of Proposition 
\ref{prop:minimal-subset} is equal to all of $\{1,\dots,k\}$ if and 
only if $\bbeta$ satisfies the stronger
inequalities of Theorem \ref{thm:main}. This requires that all $\beta_i$
are strictly positive, and if we restrict our attention to $\bbeta$ with all
$\beta_i>0$, then the condition that $J=\{1,\dots,k\}$ is 
equivalent in the polymatroid language to saying that $\bbeta$
is a ``circuit.'' 
Thus, among the vectors $\bbeta$ with $|\bbeta|=r+1$ and all $\beta_i$
strictly positive, the set of vectors satisfying
the inequalities of Theorem \ref{thm:main} is exactly the set of
circuits of the polymatroid determined by the multidegree of $X$.
% stronger inequality also implies that $\beta_i \leq n_i$ for all $i$,
% so no issue this case.
See Figure \ref{fig} for examples, and \S 1.2 of \cite{m-p-s1} for more on 
the polymatroid terminology.
\begin{figure}
\begin{center}

\begin{tikzpicture}[scale=.8]
\draw [thick, <->] (6,0)  -- (0,0) -- (0,5);
\node [below] at (6,0) {$\beta_1$};
\node [left] at (0,5) {$\beta_2$};
\draw[help lines] (0,0) grid (6,5);

\draw[fill opacity=0.8, fill=orange] (0,0) -- (0,3) -- (2,3) -- (4,1) -- (4,0)
-- cycle;
\draw [fill = gray] (0,0) circle [radius=.03];
\draw [fill = gray] (1,0) circle [radius=.03];
\draw [fill = gray] (2,0) circle [radius=.03];
\draw [fill = gray] (3,0) circle [radius=.03];
\draw [fill = gray] (4,0) circle [radius=.03];
\draw [fill = gray] (0,1) circle [radius=.03];
\draw [fill = gray] (1,1) circle [radius=.03];
\draw [fill = gray] (2,1) circle [radius=.03];
\draw [fill = gray] (3,1) circle [radius=.03];
\draw [fill = gray] (4,1) circle [radius=.08]; % base
\draw [fill = gray] (0,2) circle [radius=.03];
\draw [fill = gray] (1,2) circle [radius=.03];
\draw [fill = gray] (2,2) circle [radius=.03];
\draw [fill = gray] (3,2) circle [radius=.08]; % base
\draw [fill = gray] (0,3) circle [radius=.03];
\draw [fill = gray] (1,3) circle [radius=.03];
\draw [fill = gray] (2,3) circle [radius=.08]; % base

\draw (2,4) -- (5,1);
\draw [fill=red!40!brown] (2,4) circle [radius=.08];
\draw [fill=red!40!brown] (5,1) circle [radius=.08];
\draw [fill = green!70!black] (3,3) circle [radius=.08];
\draw [fill = green!70!black] (4,2) circle [radius=.08];
\end{tikzpicture}
\begin{tikzpicture}[scale=.5, line join=bevel, x={(.8cm, -.8cm)},z={(1.3cm,
.4cm)}]%\draw[help lines] (0,0) grid (6,5);

\draw [thick, ->] (0,0,0) -- (5,0,0);
\draw [thick, ->] (0,0,0) -- (0,5,0);
\draw [thick, ->] (0,0,0) -- (0,0,6);

\node [below] at (5,0,0) {$\beta_1$};
\node [left] at (0,5,0) {$\beta_2$};
\node [right] at (0,0,6) {$\beta_3$};

\draw[fill opacity=0.0, fill=orange] (0,0,0) -- (0,0,3) -- (2,0,3) -- (3,0,2)
-- (3,0,0) -- cycle;
\draw[fill opacity=0.0, fill=orange] (0,0,0) -- (0,3,0) -- (0,3,2) -- (0,2,3)
-- (0,0,3) -- cycle;
\draw[fill opacity=.9, fill=orange] (3,2,0) -- (2,3,0) -- (2,3,1) -- (3,2,1)
-- cycle;
\draw[fill opacity=.9, fill=orange] (3,0,0) -- (3,2,0) -- (3,2,1) -- (3,1,2)
-- (3,0,2) -- cycle;
\draw[fill opacity=.9, fill=orange] (0,3,0) -- (0,3,2) -- (1,3,2) -- (2,3,1)
-- (2,3,0) -- cycle;
\draw[fill opacity=.0, fill=orange] (0,0,3) -- (2,0,3) -- (2,1,3) -- (1,2,3)
-- (0,2,3) -- cycle;
\draw[fill opacity=0.9, fill=orange] (2,0,3) -- (3,0,2) -- (3,1,2) -- (2,1,3)
-- cycle;
\draw[fill opacity=0.9, fill=orange] (0,3,2) -- (0,2,3) -- (1,2,3) -- (1,3,2)
-- cycle;
\draw[fill opacity=.9, fill=orange] (0,0,0) -- (3,0,0) -- (3,2,0) -- (2,3,0)
-- (0,3,0) -- cycle;
\draw[fill opacity=0.9, fill=orange!60!gray] (1,3,2) -- (2,3,1) -- (3,2,1) -- (3
,1,2)
-- (2,1,3) -- (1,2,3) -- cycle; % base of polymatroid

\draw[fill opacity=0.8, fill=red!40!brown] (1,2,4) -- (1,4,2) -- (2,4,1) --
(4,2,1) -- (4,1,2) -- (2,1,4)  -- cycle;

\draw[fill opacity=0.8, fill=green!70!black] (2,2,3) -- (2,3,2) -- (3,2,2) --
cycle;

\draw [ultra thin, dashed] (1,2,4) -- (1,2,3);
\draw [ultra thin, dashed] (1,4,2) -- (1,3,2);
\draw [ultra thin, dashed] (2,4,1) -- (2,3,1);
\draw [ultra thin, dashed] (4,2,1) -- (3,2,1);
\draw [ultra thin, dashed] (4,1,2) -- (3,1,2);
\draw [ultra thin, dashed] (2,1,4) -- (2,1,3);

\end{tikzpicture}
\end{center}
\caption{Two polymatroids. The sets of bases (corresponding to our 
multidegree supports) are in gray; while the sets of circuits and
of non-circuit $1$-deficient vectors (both satisfying $|\bbeta|=r+1$) are in 
green and red, respectively.}
% first figure has delta(1)=4, delta(2)=3, delta(1,2)=5.
% 1-deficient this case (w/sum 6) automatically implies positivity
% second fig bit hard to tell, but seems similar situation.
\label{fig}
\end{figure}
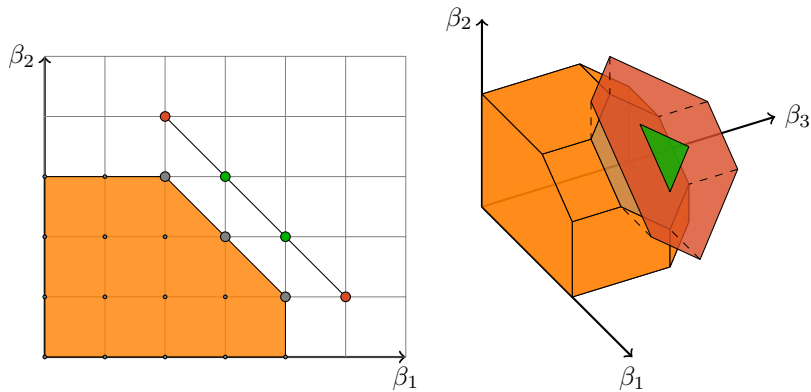
\end{rem}

\section{Slicing by products of linear spaces}

In this section, we carry out our fundamental analysis of the behavior
of slicing with products of general linear spaces. We begin with the
following.

\begin{prop}\label{prop:hypersurface} Let $X \subseteq \PP^{n_1} \times \cdots \times \PP^{n_k}$ be a 
  projective variety of dimension $r$, and suppose we have
  $\bbeta=(\beta_1,\dots,\beta_k)$ with $0 \leq \beta_i \leq n_i$ for
  $i=1,\dots,k$ and $|\bbeta| = r+1$. Write 
  $\alpha_i=n_i-\beta_i$ for each $i$.

  Consider the closed subset
$$Z_{X,\bbeta} 
=\{(L_1,\dots,L_k):X \cap (L_1 \times \cdots \times L_k)\neq \emptyset \} 
\subseteq 
{\rm Gr}(\alpha_1, n_1) \times \cdots \times {\rm Gr}(\alpha_k, n_k).$$ 
Then $Z_{X,\bbeta}$ is a hypersurface if and only if for every 
nonempty $I \subseteq \{1,\dots, k\}$ we have 
\begin{equation}\label{eq:hypersurface-ineq}
\dim p_I(X) \geq |\bbeta_I|-1,
\end{equation} 
where $p_I(X)$ denotes the 
  projection of $X$ onto $\prod_{i \in I} \PP^{n_i}$.
\end{prop}

\begin{rem}\label{rm:bounds} \rm Since $|\bbeta| = r+1$, we have 
that \eqref{eq:hypersurface-ineq} is
equivalent to having $r-\dim p_I(X) \leq |\bbeta_{I^c}|$, where 
$I^c=\{1,\dots, k\}\setminus I$. Hence, when the conditions from the 
previous Proposition are satisfied, we have that 
$r-\dim p_{I^c}(X) \leq |\bbeta_I| \leq \dim p_I(X)+1$.
The former inequality has the geometric interpretation that the
generic fiber of $X$ under $p_{I^c}$ has dimension at most $|\bbeta_I|$.
\end{rem}

\begin{proof}
  Define the incidence correspondence $V_X \subseteq X \times
{\rm Gr}(\alpha_1, n_1) \times \cdots \times {\rm Gr}(\alpha_k, n_k)$
given by
\begin{equation}\label{eq:incidence}
V_X=\{(P,L_1,\dots,L_k): P \in L_1 \times \cdots \times L_k\}
\end{equation} 
  Thus, $Z_{X,\bbeta}$ is the image of $V_X$ under projection to the
product of Grassmannians.
  Considering the projection of $V_X$ to $X$, we see that $V_X$
  is irreducible of dimension
  $d_1 = r + \sum_{i=1}^k \alpha_i(n_i-\alpha_i)$. % condition of containing
% point reduces subspace and ambient dim both by 1.
  In particular, $Z_{X,\bbeta}$ is automatically irreducible.
  The dimension of the product of Grassmannians is
  $d_2 =\sum_{i=1}^k (\alpha_1+1)(n_1-\alpha_1)$.  Because
  $d_2-d_1=1$, we will have that $Z_{X,\bbeta}$ is a hypersurface 
if and only if $V_X$ has generic fiber dimension $0$ under the 
projection to the product of Grassmannians, or equivalently, if there 
exist $L_1,\dots,L_k$ such that
  $X \cap (L_1 \times \cdots \times L_k)$ is finite and nonempty.

  First suppose that we have $\dim p_I(X) < |\bbeta_I|-1$
  for some $I \subseteq \{1,\dots, k\}$.
  Let $L_1,\dots,L_k$ be such that
  $X \cap (L_1\times \cdots \times L_k)$ is not empty, and write
$L_I =\prod _{i \in I} L_i$, and similarly for $L_{I^c}$, so that 
  $X \cap (L_1\times \cdots \times L_k) 
= (X \cap p_I^{-1} L_I) \cap p_{I^c}^{-1} L_{I^c}$.
Since
  $X \cap p_I^{-1} L_I$ is not empty, its dimension is at least the 
dimension of the generic fiber of $X$ under $p_I$, that is $r-\dim p_I(X)$.
  Hence, 
$\dim(X \cap p_I^{-1} L_I) \geq r-\dim p_I(X)> r+1-|\bbeta_I|=
|\bbeta_{I^c}|$, while
  $\codim(p_{I^c}^{-1} L_{I^c}) = |\bbeta_{I^c}|$. 
We conclude that
  $X \cap (L_1 \times \cdots \times L_k)$ has positive dimension, 
  and thus that all the inequalities are necessary in order for 
$Z_{X,\bbeta}$ to be a hypersurface.
  
Conversely, suppose that the stated inequalities are satisfied, and fix 
a point $P=(P_1,\dots,P_k) \in X$ such that 
for all $I \subsetneq \{1,\dots,k\}$, the fiber $p_I^{-1}(p_I(P)) \cap X$
has dimension less than or equal to $\sum_{i \not \in I} (n_i-\alpha_i)$.
Note that this is always possible since our inequalities are equivalent
to assuming that the minimal fiber dimension $r-\dim p_I(X)$ is less than
or equal to the desired value. Then we claim that a general
choice of $L_1,\dots,L_k$ with $P_i \in L_i$ for each $i$ will have
$X \cap (L_1 \times \cdots \times L_k)$ nonempty and finite. We prove this
by considering the incidence correspondence 
$Y \subseteq X \times {\rm Gr}(\alpha_1,n_1) \times \cdots \times {\rm Gr}(\alpha_k,n_k)$
consisting of $(Q,L_1,\dots,L_k)$ with $L_i \ni P_i$ for all $i$, and
$Q \in L_1 \times \cdots \times L_k$. Considering the case $Q=P$, we
see that the image of $Y$ under projection to 
${\rm Gr}(\alpha_1,n_1) \times \cdots \times {\rm Gr}(\alpha_k,n_k)$ is exactly the tuples
of $L_i$ containing $P_i$, which is itself a product of Grassmannians,
having dimension $\sum_i \alpha_i(n_i-\alpha_i)$. The finiteness statement
we want amounts to showing that $Y$ has generically finite fibers under 
this projection,
or equivalently, that the dimension of $Y$ is no larger than the dimension
of its image. We do this by decomposing $Y$ into locally closed subsets
$Y_I$, defined to be the subset of $Y$ on which $Q_i=P_i$ precisely
when $i \in I$ (here we allow $I=\emptyset$). We then consider the 
projection of $Y_I$ onto $X$. First,
$Y_I \neq \emptyset$ if and only if $\alpha_i\geq 1$ for all $i \not \in I$.
In this case, $Y_I$ maps into $p_I^{-1}(p_I(P)) \cap X \subseteq X$,
% actually surjects onto explicit open subset
and every fiber will have dimension equal to
$$\sum_{i \in I} \alpha_i(n_i-\alpha_i)+
\sum_{i \not \in I} (\alpha_i-1)(n_i-\alpha_i).$$
We conclude that $Y_I$ has dimension less than or equal to 
$$\sum_{i \not \in I}(n_i-\alpha_i) +\sum_{i \in I} \alpha_i(n_i-\alpha_i)+
\sum_{i \not \in I} (\alpha_i-1)(n_i-\alpha_i)=
\sum_i \alpha_i(n_i-\alpha_i),$$
as desired.
\end{proof}

\begin{ex} Consider a product of varieties 
$X = X_1 \times \cdots \times X_k$, so that each $X_i \subseteq \PP^{n_i}$ 
has dimension $r_i$, and $\sum_i r_i = r$. If $\bbeta$ is 
such that $|\bbeta|= r+1$, and $Z_{X,\bbeta} \subseteq 
{\rm Gr}(\alpha_1, n_1) \times \cdots \times {\rm Gr}(\alpha_k, n_k)$ is 
as defined above, then Proposition~\ref{prop:hypersurface} states that 
$Z_{X,\bbeta}$ has codimension $1$ if and only if 
$\sum_{i \in I} r_i \leq |\bbeta_I| \leq \sum_{i \in I} r_i +1$ 
for all $I$. This can occur if and only if $\beta_j = r_j + 1$ for one 
index and $\beta_i = r_i$ otherwise (a condition that can also be deduced 
directly from $Z_{X,\bbeta}$). In this case, 
$Z_{X,\bbeta}=
{\rm Gr}(\alpha_1,n_1)\times \cdots \times Z_{X_j} \times \cdots 
\times {\rm Gr}(\alpha_k, n_k)$, where $Z_{X_j}$ is the hypersurface arising
in the classical Cayley-Chow construction.
\end{ex}

\begin{ex} Consider a variety $X_0 \subseteq \PP^{n_0}$, of dimension $r_0$,
and let $X$ be the image of $X_0$ in the diagonal embedding 
$\PP^{n_0} \rightarrow (\PP^{n_0})^k$. In this case, the conditions in 
Proposition~\ref{prop:hypersurface} require that 
$0 \leq |\bbeta_I| \leq r_0 + 1$. Hence, 
$Z_{X,\bbeta}$ is a hypersurface for any choice of 
$(\beta_1,\dots,\beta_k)$ summing to $r_0 + 1$. 
More specifically, we have a surjective rational map
$${\rm Gr}(\alpha_1,n_1) \times \cdots \times {\rm Gr}(\alpha_k,n_k)
\dashrightarrow {\rm Gr}(n_0-r_0-1,n_0)$$
given by intersection of linear spaces, and it is clear that on the open
subset where this map is defined, we have that $Z_{X,\bbeta}$ is 
the preimage of the hypersurface $Z_{X_0}$ arising in the classical 
Cayley-Chow construction. By irreducibility of $Z_{X,\bbeta}$,
we see that it must simply be the closure of the preimage of $Z_{X_0}$.
\end{ex}

\begin{ex}\label{ex:decomposable} Suppose that for some 
$I \subsetneq \{1,\dots,k\}$, we have $X' \subseteq \prod_{i \in I} \PP^{n_i}$
such that $X=X' \times \prod_{i \not \in I} \PP^{n_i}$ (equivalently,
$X=p_I^{-1}(p_I(X))$). Write $r'$ for the dimension of $X'$, and
suppose we have $\bbeta$ satisfying \eqref{eq:hypersurface-ineq}.
Then we have $|\bbeta_I| \leq \dim p_I(X)+1=r'+1$, so we must
have $|\beta_{I^c}| \geq \sum_{i \in I^c} n_i$, and then it follows
that $\beta_i=n_i$ for all $i \in I^c$, and also that $|\beta_I|=r'+1$.
In particular, the inequalities of Theorem \ref{thm:main} are violated.
However, one easily verifies that 
$$Z_{X,\bbeta}=Z_{X',\bbeta_I} \times 
\prod_{i \not\in I} {\rm Gr}(\alpha_i, n_i),$$
so the study of $Z_{X,\bbeta}$ in this case reduces to the study of
$Z_{X',\bbeta_I}$, and in particular $X$ can be recovered from $Z_{X,\bbeta}$
if and only if $X'$ can be recovered from $Z_{X',\bbeta_I}$.
% X of this form iff \dim p_I(X)=\dim X - \sum_{i \in I^c} n_i, so
% again this read off from multidegree.
\end{ex}

We now analyze when $X$ can be recovered from $Z_{X,\bbeta}$.

\begin{prop}\label{prop:degenerate} In the situation of 
Proposition~\ref{prop:hypersurface} (and in particular assuming
\eqref{eq:hypersurface-ineq}), 
we have that $X$ is uniquely determined by $Z_{X,\bbeta}$ if 
\begin{equation}\label{eq:good-ineqs}
\dim p_I(X) \geq |\bbeta_I|
\end{equation}
for all $I \subsetneq \{1,\dots, k\}$.

Conversely, if $X$ is not of the form of Example \ref{ex:decomposable},
and $X$ is uniquely determined by $Z_{X,\bbeta}$, then \eqref{eq:good-ineqs}
is satisfied. 
\end{prop} 

\begin{proof} 
Let $S_Z$ be the set of points $(P_1,\dots, P_k)$ with the property that
every $L_1 \times  \cdots \times L_k$ containing $(P_1,\dots, P_k)$ has
$X \cap (L_1 \times \cdots \times L_k) \neq \emptyset$. Then obviously
$X \subseteq S_Z$. We claim that if \eqref{eq:good-ineqs} is satisfied, then
\begin{equation}
\label{eq:degenerate}
S_Z \subseteq X \cup \bigcup_{I \subsetneq \{1,\dots, k\}} p_I^{-1}(X_I)
\end{equation}
where $X_I \subseteq p_I(X)$ is the closed subset over which the 
fibers of $X$ under $p_I$ have dimension greater than or equal to
$\sum_{i \not \in I} (n_i-\alpha_i)$. Note that \eqref{eq:good-ineqs}
implies that $X_I \neq p_I(X)$.
From the claim we see that -- given the multidegree of $X$ --
we can recover $X$ from $S_Z$: indeed, by Corollary \ref{cor:multdeg-projs}
the multidegree of $X$ determines $\dim p_I(X)$ for all $I$, and we see
that every other potential component of $S_Z$ has dimension strictly smaller
than $X$ under at least one projection.

To prove the claim, fix $P=(P_1,\dots,P_k) \not \in X$ and with 
$p_I(P) \not \in X_I$ for all $I \subsetneq \{1,\dots,k\}$. We wish to
show that $P \not \in S_Z$, or equivalently, that there exist
$L_1,\dots,L_k$ with $P_i \in L_i$ for all $i$, and with
$X \cap (L_1 \times \cdots \times L_k) = \emptyset$. The proof is 
similar to the proof of Proposition \ref{prop:hypersurface}:
consider the incidence correspondence 
$Y \subseteq X \times {\rm Gr}(\alpha_1,n_1) \times \cdots \times {\rm Gr}(\alpha_k,n_k)$
consisting of $(Q,L_1,\dots,L_k)$ with $L_i \ni P_i$ for all $i$, and
$Q \in L_1 \times \cdots \times L_k$. In this case, we wish to show
that the image of $Y$ under projection to 
${\rm Gr}(\alpha_1,n_1) \times \cdots \times {\rm Gr}(\alpha_k,n_k)$ does not contain all
tuples of $L_i$ containing $P_i$, so it will suffice to show that $Y$
has dimension strictly smaller than $\sum_i \alpha_i(n_i-\alpha_i)$.
We decompose $Y$ into the subsets $Y_I$ as before. For 
$I\subsetneq \{1,\dots,k\}$, just as before
$Y_I$ has image contained in $p_I^{-1}(p_I(P)) \cap X \subseteq X$,
with every fiber having dimension equal to
$$\sum_{i \in I} \alpha_i(n_i-\alpha_i)+
\sum_{i \not \in I} (\alpha_i-1)(n_i-\alpha_i).$$
We conclude that $Y_I$ has dimension less than or equal to 
$$\sum_{i \not \in I}(n_i-\alpha_i)-1 +\sum_{i \in I} \alpha_i(n_i-\alpha_i)+
\sum_{i \not \in I} (\alpha_i-1)(n_i-\alpha_i)=
\sum_i \alpha_i(n_i-\alpha_i)-1,$$
as desired. 

Conversely, suppose that we have some $I \subsetneq \{1,\dots,k\}$ 
such that $|\bbeta_I| > \dim p_I(X)$. %; then necessarily 
%$\sum_{i \in I} \beta_i = \dim p_I(X)+1$. 
Then we have that for a
general choice of $L_i$ for $i \in I$, the intersection 
$X \cap \bigcap_{i \in I}p_i^{-1} L_i$ is empty.
% could cite hypersurface prop, but kind of overkill: every time cut by
% pullback of general hyperplane, dimension has to go down by 1
Thus, $p_I(Z_{X,\bbeta})$ is a proper subset of 
$\prod_{i \in I} {\rm Gr}(\alpha_i, n_i)$. Since $Z_{X,\bbeta}$ is 
a hypersurface, it must be equal to 
$p_I(Z_{X,\bbeta}) \times 
\prod_{i \not\in I} {\rm Gr}(\alpha_i, n_i)$. We conclude that $Z_{X,\bbeta}$ 
is invariant under applying automorphisms of $\PP^{n_i}$ for $i \not\in I$, 
and since we have assumed that $X$ is not of the form of Example
\ref{ex:decomposable}, it follows that $X$ cannot be recovered from 
$Z_{X,\bbeta}$.
\end{proof}

Note that (unlike in the classical case), it may really be the case
that the $S_Z$ in the above proof contains components other than $X$;
see Example \ref{ex:trifocal} below.

\begin{rem} In the case that the inequalities \eqref{eq:hypersurface-ineq}
are satisfied, we can understand the hypersurface $Z_{X,\bbeta}$ 
as follows:
according to Proposition \ref{prop:minimal-subset}, there a nonempty
$J \subseteq \{1,\dots,k\}$ which is minimal -- in the strong sense --
satisfying 
$$|\bbeta_J| = \dim p_J(X)+1.$$ 
In this case, $X'=p_J(X)$ will satisfy our setup with the stronger 
inequalities \eqref{eq:good-ineqs}, so $X'$ can be recovered from
$Z_{X',\bbeta_J}$. Furthermore, we will have $Z_{X,\bbeta}$ 
equal to the product of $Z_{X',\bbeta_J}$ with all the 
${\rm Gr}(\alpha_i,n_i)$ for $i \not \in J$, so the information in 
$Z_{X,\bbeta}$ is exactly equal to the data of $p_J(X)$.
\end{rem}

\section{Tensor products of unique factorization domains}

In the classical setting, the fact that the hypersurface $Z_X$ is the
zero set of a single polynomial $F_X$ in Pl\"ucker coordinates is a 
consequence of the fact that the homogeneous coordinate ring of a 
Grassmannian in Pl\"ucker coordinates is a unique factorization domain (UFD).
We will make use of this to conclude the same statement in our case, 
but this requires a certain amount of care, as the condition of being
a UFD is not very stable (for instance, there
are examples where $R$ is a UFD, but the power
series right $R[[t]]$ is not). We address this with the following
proposition, which states that under relatively mild additional
hypotheses, Gauss' argument for unique factorization in a polynomial
ring over a UFD extends to more general tensor products.

For the following proposition, we temporarily drop the hypothesis that
we are working over a field $K$.

\begin{prop}\label{prop:ufd-tensor}
Let $A$ be a ring, $B$ and $C$ algebras over $A$,
and suppose that $B$ is a Noetherian UFD, $C$ is flat and finitely
generated over $A$, and for every field 
$K$ over $A$, we have that $K \otimes_A C$ is a UFD, with unit group equal 
to $K^{\times}$. Then $B \otimes_A C$ is a UFD.
\end{prop}

\begin{proof} 
First observe that under our hypotheses, if $B'$ is an $A$-algebra with
fraction field $K'$, then we have injections 
$K' \to K' \otimes_A C$ and $B' \to B' \otimes_A C \to K' \otimes_A C$. 
Indeed, injectivity of the last map follows from
flatness of $C$, while injectivity of the first is a consequence of
the implicit hypothesis that $K' \otimes_A C$, being a domain, is not
the zero ring. % Thus, the first map is not zero, and must be injective.
Injectivity of the map $B' \to B' \otimes_A C$ follows from the 
injectivity of the first map. Note also that $B \otimes_A C$ is
finitely generated over a Noetherian ring, hence Noetherian, so it
suffices to show that every irreducible element is prime.

Our first claim is that if $x \in B$ is a prime element, then
$x \otimes 1$ is prime in $B \otimes_A C$, and given $y \in B$ we have
that $y$ is a multiple of $x$ if and only if $y \otimes 1$ is a multiple
of $x \otimes 1$ in $B \otimes_A C$. Indeed, if we let $B'=B/(x)$ and
apply the above, the second statement follows immediately
from the injectivity of 
$B/(x) \to B/(x) \otimes_A C = (B \otimes_A C)/(x \otimes 1)$,
while the injection $B/(x) \otimes_A C \to K' \otimes_A C$ together
with the hypothesis that $K' \otimes_A C$ is an integral domain implies
that $x \otimes 1$ is prime. From the second part of the claim, we can
conclude that if $K$ is the fraction field of $B$, then the intersection of 
$K$ with $B \otimes_A C$ inside of $K \otimes_A C$ is equal to $B$.
% just clear denominators and apply to primes in denominator
We also see that conversely if $x \otimes 1$ is irreducible in $B \otimes_A C$
then $x$ must be irreducible in $B$.
%otherwise could write as product involving at least two prime factors,
% which stay prime in tensor product.
In this case, $x$ is prime in $B$, and hence $x \otimes 1$ is prime in
$B \otimes_A C$.

Now, suppose $x$ is irreducible in $B \otimes_A C$,
and consider the image of $x$ in $K \otimes_A C$.
If $x$ becomes a unit in $K \otimes_A C$, then by hypothesis it is of the 
form $y \otimes 1$ for some $y \in K^{\times}$, so we have from the above
that $y \in B$, so $x = y \otimes 1$ is prime in $B \otimes_A C$.

Next, we see that any nonzero element of $K \otimes_A C$ can be written 
uniquely up to $B^\times$ in the form $\alpha f$, for $\alpha \in K^*$, 
and $f \in B \otimes_A C$
with the property that $f$ is not a multiple of any non-unit in $B$.
Indeed, every element of $K \otimes_A C$ can be multiplied by an element
of $B$ to clear denominators,
% since any tensor is a finite sum
so is of the form $\alpha f$ where $f \in B \otimes_A C$. Obviously, if
$f$ is a multiple of a non-unit in $B$, we can absorb it into $\alpha$,
% need process terminates -- but needed that anyhow. know non-unit in B 
% is nonunit in tensor prod since showed primes stay prime
so we can write the element in the desired form. Uniqueness follows
by taking two such representations, clearing denominators, and using that
primes in $B$ remain prime in $B \otimes_A C$. 

It then follows that if we have an irreducible element of $B \otimes_A C$ 
which does not become
a unit in $K \otimes_A C$, then it must remain irreducible in 
$K \otimes_A C$. Indeed, a nontrivial factorization in $K \otimes_A C$ can be
represented as $(\alpha_1 f_1)(\alpha_2 f_2)$ as above, and then we
see that since $f_1$ and $f_2$ are not multiples of any non-units of $B$,
the same is true of $f_1 f_2$. 
% again using primes in B remain prime.
Thus, the hypothesis that $\alpha_1 \alpha_2 f_1 f_2 \in B \otimes_A C$
implies that $\alpha_1 \alpha_2$ can't have any denominators, and
then $(\alpha_1 \alpha_2 f_1) f_2$ gives a nontrivial factorization in
$B \otimes_A C$. Finally, we conclude that any such irreducible element must 
be prime: it is prime in $K \otimes_A C$ by hypothesis, so if it divides
a product in $B \otimes_A C$, it divides one of the factors in 
$K \otimes_A C$. But again using the above representation, we see that
it must also divide the same factor in $B \otimes_A C$.
% can write factor as original times $\alpha f$, again conclude $\alpha$
% can't have denominators since original has no $\alpha$ term.
We thus conclude that every irreducible element of $B \otimes_A C$ is
prime, and hence that $B \otimes_A C$ is a UFD.
\end{proof}

Returning to varieties over $K$, we then conclude the desired statement 
on hypersurfaces in products of Grassmannians.

\begin{cor}\label{cor:multiplucker} Let 
$G:=G_1 \times \cdots \times G_k \subseteq \PP^{N_1} \times \cdots \times
\PP^{N_k}$
be a product of Pl\"ucker embeddings of Grassmannians, and let
$Z \subseteq G$ be a hypersurface. Then
$Z=Z(F)$ for some multihomogeneous form $F$.
\end{cor}

\begin{proof} First, we have that $Z$ corresponds to a multihomogeneous
prime ideal of height $1$ in $S(G)$, the multihomogeneous coordinate ring
of $G$. 
% $I(Z)$ is multihomogeneous by defn. can prove height 1 by looking at
% cones as in exers 6.1.9 & 6.1.10 of vars notes. looking at ambient spaces,
% have morphism between nonsingular vars w/constant-dimensional irred fibers,
% so flat, open, get irreds correspond. rest of arg goes through as before.
% preimages not so related to undefined locus, but fiber dim is k, so works
% out same.
Now, we claim that $S(G)$ is a UFD. 
Since $S(G)=S(G_1) \otimes_K \cdots \otimes_K S(G_k)$, we will prove
this by induction on $k$, using Proposition \ref{prop:ufd-tensor}.
The base case is exactly the classical case; see Proposition 2.1 of
Chapter 3 of \cite{g-k-z1}. 
Thus, we need only observe that $S(G_k)$ over $K$ satisfies the 
hypotheses of the proposition, most of which are immediate: flatness
is automatic over $K$, the hypothesis on the units comes from the fact
that $S(G_k)$ is the homogeneous coordinate ring of a projective variety,
and the hypothesis that for any field extension $K'$ over $K$ we have
that $K' \otimes_K S(G_k)$ is a UFD also follows from
the classical case, since $K' \otimes_K S(G_k)$ is simply the homogeneous
coordinate ring of the Pl\"ucker embedding over $K'$.
% Noetherian/finite generation hypotheses are obvious in context
We thus conclude that $S(G)$ is a UFD, and therefore that $Z=Z(F)$ for some 
$F \in S(G)$.  Finally, $F$ must be multihomogeneous, or it could not
generate a multihomogeneous prime ideal.
% eg, differs from multihomogeneous generator by element of $K^{\times}$.
\end{proof}

\section{(Multi)degrees and Cayley-Chow forms}

We are now ready to define multigraded Cayley-Chow forms. In order to 
obtain good behavior of multidegrees, there is one additional twist
to consider. Namely, unlike in the classical case, it is possible that
the Cayley-Chow form naturally has a multiplicity greater than $1$.

\begin{defn}\label{def:multichow} In the situation of Proposition 
\ref{prop:hypersurface}, suppose also that \eqref{eq:hypersurface-ineq}
is satisfied. Let $\epsilon_{X,\bbeta}$ be the degree of the
map $V_X \to Z_{X,\bbeta}$, where $V_X$ is as in \eqref{eq:incidence}.

Then let $F_X$ be a multihomogeneous polynomial in multi-Pl\"ucker coordinates
with $Z(F_X)=Z_{X,\bbeta}$ (Corollary \ref{cor:multiplucker}), 
and define the \textbf{multigraded Cayley-Chow form} of $X$ to be 
$$H_{X,\bbeta}:= F_X^{\epsilon_{X,\bbeta}}.$$
\end{defn}

Recall that $V_X$ is the incidence correspondence consisting of a point
of $X$ together with a tuple of linear spaces containing the coordinates
of the point. Then the map to $Z_{X,\bbeta}$ is simply the map forgetting
the point of $X$.

Without further hypotheses, it may certainly be the case that 
$\epsilon_{X,\bbeta}>1$. 

\begin{ex}\label{ex:prod-curves} Let 
$X=C_1 \times C_2 \subseteq \PP^2 \times \PP^2$,
where each $C_i$ is the curve defined by a homogeneous form $F_i$ of degree 
$d_i$. In order for
\eqref{eq:hypersurface-ineq} to be satisfied, we need to have either
$\bbeta=(1,2)$ or $\bbeta=(2,1)$. In the first case,
we see that if we have fixed $L,P$ with 
$X \cap (L \times P) \neq \emptyset$, then in fact $X \cap (L \times P)$
contains $d_1$ points, so $\epsilon_{X,\bbeta}=d_1$. Meanwhile
$Z_{X,\bbeta}$ depends only on $P$, and the $F_X$ of Definition 
\ref{def:multichow}
is simply $F_2$. Thus, $H_{X,\bbeta}=F_2^{d_1}$. Similarly,
if $\bbeta=(2,1)$, we find 
$H_{X,\bbeta}=F_1^{d_2}$.
\end{ex}

For a more interesting example in positive characteristic, see
Example \ref{ex:frobenius} below. However, we have the following.

\begin{prop}\label{prop:degree-1} In the situation of
Definition \ref{def:multichow}, suppose further that the inequalities of
\eqref{eq:good-ineqs} are satisfied. Then the map 
$V_X \to Z_{X,\bbeta}$ is
generically injective, and if $K$ has characteristic $0$, we have
$\epsilon_{X,\bbeta}=1$.
\end{prop}

\begin{proof} The generic injectivity amounts to saying that if
$P=(P_1,\dots,P_k) \in X$ is general, then general choices of $L_i$ containing
$P_i$ will have $X \cap (L_1 \times \cdots \times L_k)=\{P\}$.
We recall that \eqref{eq:good-ineqs} implies that for all 
$I \subsetneq \{1,\dots,k\}$, the generic fiber dimension of $X$ under
$p_I$ is strictly less than $|\bbeta_{I^c}|$. We prove by
induction on $k$ the following slightly more general statement: suppose that 
$Y \subseteq \prod_i \PP^{n_i}$ is a pure-dimensional algebraic set such 
that for every $I \subsetneq \{1,\dots,k\}$, every component of $Y$ has 
generic fiber dimension under $p_I$ strictly less than
$|\bbeta_{I^c}|$. Then there exists a dense open subset $U$
of $Y$ such that for every $P=(P_1,\dots,P_k) \in U$, a general choice of
$L_i$ containing $P_i$ will have 
$Y \cap (L_1 \times \cdots \times L_k)=\{P\}$.
Note that we allow $I=\emptyset$ in our hypotheses, which says simply 
that $Y$ has dimension strictly smaller than $|\bbeta|$. 
% this version bit confusing: general fiber may not be pure in reducible
% case, but since reduce to irreducible case before passing to fiber, seems OK
% version w/stronger hyps below:
% $Y \subseteq \prod_i \PP^{n_i}$ is an algebraic set such that for every 
% $I \subsetneq \{1,\dots,k\}$, every component of $Y$ has the same generic
% fiber dimension under $p_I$, and this dimension is strictly less than
% $|\bbeta_{I^c}|$. Then there exists a dense open subset $U$
% of $Y$ such that for every $P=(P_1,\dots,P_k) \in U$, a general choice of
% $L_i$ containing $P_i$ will have 
% $Y \cap (L_1 \times \cdots \times L_k)=\{P\}$.
% Note that we allow $I=\emptyset$ in our hypotheses, which says simply 
% that $Y$ is pure-dimensional, of dimension strictly smaller than 
% $|\bbeta|$. 
% NB: without purity hypotheses, this fails: eg take Y to be some X together
% with a point in S_Z. Even sats comp-by-comp conds on dims of generic fibers

We first observe that the desired statement reduces to the case that 
$Y$ is irreducible. Indeed, if $Y_1,\dots,Y_n$ are the components of $Y$, 
and if we construct $S_{Z,i}$ from each $Y_i$ as in the proof of 
Proposition \ref{prop:degenerate} (with $Y_i$ in place of $X$), then
\eqref{eq:degenerate} together with our hypotheses on the $Y_i$ implies
that we cannot have $Y_i \subseteq S_{Z,i}$ for any distinct $i, j$.
To see this, if we write $Y_{I,i}$ in place of $X_I$, and write 
$r=\dim Y=\dim Y_i$, we see that each
$Y_{I,i}$ must have dimension strictly less than $r-|\bbeta_{I^c}|$,
% since dimension of preimage is at least $\dim Y_{I,i}+|\bbeta_{I^c}|$
% by definition, but must be proper in $Y_i$.
while our hypotheses imply that $\dim p_I(Y_j)>r-|\bbeta_{I^c}|$, so
we cannot have $p_I(Y_j) \subseteq Y_{I,i}$.
Thus, if we suppose that $U_i$ satisfies the desired conditions for each $Y_i$
separately, we then see that 
$$\bigcup_i \left(U_i \smallsetminus \bigcup_{j \neq i} S_{Z,j}\right)$$
satisfies the desired condition for all of $Y$.
Consequently, for simplicity we will henceforth assume $Y$ is irreducible.

Then the base case our of induction is $k=1$, which is exactly the classical 
situation, and is proved by considering projection from any point in $Y$.
For induction, we claim that the general fiber of $Y$ under $p_k$
satisfies our hypotheses for $k-1$. First, since $Y$ is irreducible the 
general fiber will be pure-dimensional. 
% by semicontinuity of fiber dimension
Next, for any fixed $P_k \in p_k(Y)$,
and $I' \subsetneq \{1,\dots,k-1\}$, if $Y_{P_k}$ denotes the fiber of $Y$
over $P_k$, and if we set $I=I' \cup \{k\}$, we see
that for any $(P_i)_{i \in I'} \in p_{I'}(Y_{P_k})$, the fiber of 
$Y_{P_k}$ over $(P_i)_{i \in I'}$ is equal to the fiber of $Y$ over 
$(P_i)_{i \in I}$.  Note that $|\bbeta_{I^c}|$ is the desired
dimension bound also for $k-1$, since now we take the complement of $I'$
in $\{1,\dots,k-1\}$. By hypothesis, there is an open subset
$U_I \subseteq p_I(Y)$ such that the fiber dimension of $Y$ over any point of
$U_I$ is strictly less than $|\bbeta_{I^c}|$. Then we observe
% version for stronger hyps above:
% $U_I \subseteq p_I(Y)$ such that the fiber dimension of $Y$ over any point of
% $U_I$ is equal to $\dim Y - \dim p_I(Y)$, which is strictly less than 
% $|\bbeta_{I^c}|$ by hypothesis. Then we observe
that there is an open subset $V_I$ of $p_k(Y)$ such that $p_I^{-1}(U_I)$ is 
dense in every fiber $Y\cap p_k^{-1}(Q)$ for $Q \in V_I$: indeed, this follows 
from constructibility of 
images, together with semicontinuity of fiber dimension, since the only 
way that $p_I^{-1}(U_I)$ can fail to be dense in the fiber over some
$Q\in p_k(Y)$ is if the
dimension of $Y\smallsetminus p_I^{-1}(U_I)$ over $Q$ is strictly larger 
than the generic fiber dimension of $Y \smallsetminus p_I^{-1}(U_I)$ over 
$p_k(Y)$. Then a general $P_k$ will not only yield $Y_{P_k}$ pure-dimensional,
but will also lie in every $V_I$, so we see that every component of 
$Y_{P_k}$ necessarily satisfies the desired generic fiber dimension bound.

Now, we note that the set of $(P,L_1,\dots,L_k) \in Y \times 
{\rm Gr}(\alpha_1, n_1) \times \cdots \times {\rm Gr}(\alpha_k, n_k)$
such that $Y \cap (L_1 \times \cdots \times L_k) =\{P\}$ is
constructible: indeed, this follows from semicontinuity of fiber dimension,
properness, and constructibility of connected fibers (see Theorem 9.7.7 of
\cite{ega43}).\footnote{Here we are considering only classical points;
the correct statement for schemes involves geometrically connected fibers,
but since we work over an algebraically closed field and consider only
classical points, this amounts to the same thing.} Thus, to prove the
desired statement, it suffices to show that this set is Zariski dense
inside the locus of points satisfying $P \in L_1 \times \cdots \times L_k$.
If this were not the case, it would be contained in a proper Zariski closed
subset $Z$. We would then have that a dense open subset of $Y$ has the
property that $Z$ does not fully contain any of the fibers over that subset:
this is, there would be a dense open subset of points $P \in Y$ such that
the choices of $L_i$ containing $P_i$ and having
$Y \cap (L_1 \times \cdots \times L_k) = \{P\}$ are contained in a proper
closed subset.
% follows from openness of projection map, considering complement of $Z$

On the other hand,
by the above claim and the induction hypothesis, if we fix a general 
$P_k \in p_k(Y)$, and let $Y_{P_k}$ be the corresponding fiber, then
we know that for $(P_1,\dots,P_{k-1})$ general in $Y_{P_k}$, and
$L_i$ general containing $P_i$ for $i=1,\dots,k-1$, we have
$Y_{P_k} \cap (L_1 \times \cdots \times L_{k-1})=\{(P_1,\dots,P_{k-1})\}$.
We also know that
$$\dim p_k(Y \cap (L_1 \times \cdots \times L_{k-1} \times
\PP^{n_k})) < \beta_k,$$
so $Y \cap (L_1 \times \cdots \times L_{k-1} \times \PP^{n_k})$
will meet a general $L_k$ containing $P_k$ only in $P_k$.
Thus, we will have
$Y \cap (L_1 \times \cdots \times L_k)=\{P\}$. Given the generality of
$P_i$ and $L_i$, this proves the desired statement.

We have thus proved the generic injectivity statement in general. To prove
that $\epsilon_{X,\bbeta}=1$ in characteristic $0$, we use
the Bertini theorem given as Corollary 5 of \cite{kl2}: a general divisor
in a basepoint-free linear system is smooth at all smooth points of the 
ambient scheme. In particular, intersecting a generically
reduced scheme with the preimage of a general hyperplane in any of the
$\PP^{n_i}$ will yield another generically reduced scheme. Applying this
inductively, if we fix general $L_1,\dots,L_{k-1}$, and general $L_k'$ of 
dimension $\alpha_k+1$, then 
$$X \cap (L_1 \times \cdots \times L_{k-1} \times L_k')$$ 
will consist of a finite number of reduced points, and the same will still
be true if we further intersect with any $L_k$ of codimension $1$
in $L_k'$. Since we have already shown that such an intersection consists of
a single point, we conclude that it is in fact a single reduced point,
and $\epsilon_{X,\bbeta}=1$.
\end{proof}

\begin{ex} The following example demonstrates the delicacy of the
inductive statement proved in Proposition \ref{prop:degree-1}:
given $d_1,d_2>1$, let $S_1, S_2$ be surfaces of degrees $d_1,d_2$ 
respectively in $\PP^3$, and $C_1, C_2$ curves in $S_1, S_2$. Set
$$Y = (C_1 \times S_2) \cup (S_1 \times C_2),$$
and $\beta_1=\beta_2=2$.
Then $p_1(Y)$ and $p_2(Y)$ are both irreducible, and the ``generic fibers''
of both $p_i$ have dimension $1$ (which is strictly smaller than either
$\beta_i$), in the sense that there are open dense subsets of each $p_i(Y)$
on which the fiber dimension is $1$. Thus, $Y$ satisfies a slight variant
of the induction statement, but we see that any $L_1 \times L_2$ meeting $Y$ 
has to contain at least $\min(d_1,d_2)$ points.
\end{ex}

In order to complete the proofs of the basic properties of our generalized
Cayley-Chow construction, it will be helpful to extend it to cycles. 
The main reason for this is that even if $X$ is irreducible, its
intersections with general products of linear spaces are 
not always irreducible (Bertini theorems imply that they usually are, 
but they will not be if for instance some projection has $1$-dimensional 
image). Some preliminary notation is the following.

\begin{notn} If $X \subseteq Y$ is a pure-dimensional closed subscheme
of a smooth variety,
denote by $[X]$ the associated cycle. If $\Xi,\Xi'$ are cycles on $Y$
which meet in the expected codimension, write $\Xi \cdot \Xi'$
for the induced intersection cycle (see for instance Serre's definition
of intersection multiplicity on p.\ 427 of \cite{ha1}).
\end{notn}

Note that we do not work up to rational (or other) equivalence; the
point of introducing the notation is that it is not always true that
$[X] \cdot [X']=[X \cap X']$, even when $X$ and $X'$ intersect in the
expected dimension. 
% this is rephrasing of need for Serre: eg intersection of non-CM union
% of planes with a third plane. can also turn this into failure of naive
% projection formula, looking at fibers of a suitable projection onto a
% plane.
In our situations, we will have 
$[X] \cdot [X']=[X \cap X']$ due to generality hypotheses, but we will
have to justify this point.

\begin{defn}\label{defn:cayley-chow-cycle}
Given $n_1,\dots,n_k$, $r$ and $\bbeta=(\beta_1,\dots,\beta_k)$ 
with $|\bbeta| = r+1$ and $0 \leq \beta_i \leq n_i$ for all $i$,
linearly extend the construction 
$X \mapsto \epsilon_{X,\bbeta}[Z_{X,\bbeta}]$ to 
effective $r$-cycles on $\PP^{n_1} \times \cdots \times \PP^{n_k}$ by
using our previous construction for those components
satisfying \eqref{eq:hypersurface-ineq}, and extending by zero for any
additional components.
Extend the resulting Cayley-Chow form construction multiplicatively.
For an effective $r$-cycle $\Xi$, denote the resulting constructions by 
$Z_{\Xi,\bbeta}$ and $H_{\Xi,\bbeta}$ respectively.
\end{defn}

Note that we are incorporating the multiplicities into our new notation,
so that if $\Xi=[X]$, we have 
$Z_{\Xi,\bbeta}=\epsilon_{X,\bbeta}[Z_{X,\bbeta}]$.

\begin{rem}\label{rem:extension-compatible} Observe that we can
rephrase Definition \ref{def:multichow} as saying that we are taking
the form cutting out $\epsilon_{X,\bbeta}[Z_{X,\bbeta}]$, and the latter
is $\pi_* [V_X]$, where $\pi$ is the projection to the product of
Grassmannians. Definition
\ref{defn:cayley-chow-cycle} allows us to extend this as follows:
if $\Xi$ is an effective $r$-cycle,
we can construct an incidence correspondence cycle $V_{\Xi}$
on $\left(\prod_{i \in I^c} \PP^{n_i}\right) \times 
\left(\prod_{i \in I^c} {\rm Gr}(\alpha_i,n_i)\right)$ by linearly extending
our previous construction, and our extension by zero in Definition
\ref{defn:cayley-chow-cycle} means that we will still have
$Z_{\Xi,\bbeta}=\pi_* V_{\Xi}$.
 Indeed, according to the proof of
Proposition \ref{prop:hypersurface}, a component of $\Xi$ fails to 
satisfy \eqref{eq:hypersurface-ineq} precisely when the corresponding 
component of $V_{\Xi}$ drops dimension under $\pi$.
\end{rem}

We then have the following description of the behavior of
multigraded Cayley-Chow forms under partial evaluation.

\begin{prop}\label{prop:form-evaluate} 
Let $\Xi$ be an effective $r$-cycle as in Definition
\ref{defn:cayley-chow-cycle}, and $H_{\Xi,\bbeta}$ its
associated multigraded Cayley-Chow form. For any $I \subsetneq \{1,\dots,k\}$,
given general $(L_i)_{i \in I}$ of codimensions $\beta_i$, set
$L_I:=\prod_{i \in I} L_i$. Then the partial evaluation of
$H_{\Xi,\bbeta}$ at the $L_i$ for $i \in I$ yields 
the multigraded Cayley-Chow form associated to
$p_{I^c *}([p_I^{-1}(L_I)] \cdot \Xi)$ and $\bbeta_{I^c}$.
\end{prop}

In the above, if $p_I^{-1}(L_I)$ does not meet the support of $\Xi$, we
should interpret the associated multigraded Cayley-Chow form to be constant.
% won't happen if stronger ineqs satisfied, but will with weaker ones.

\begin{rem}\label{rem:intersect-agree} Note that if $X$ is a subvariety,
the generality of $L_I$ implies that $p_I^{-1}(L_I) \cap X$ is 
necessarily pure-dimensional, of codimension $|\bbeta_I|$ in $X$, so it is 
reasonable to pass to the associated cycle, and apply pushforward of 
cycles.
If $X$ is not Cohen-Macaulay, then even when $p_I^{-1}(L_I)$ meets $X$ in
the expected codimension, we could \emph{a priori} have 
that $[p_I^{-1}(L_I) \cap X] \neq [p_I^{-1}(L_I)] \cdot [X]$,
so we have to be slightly careful with our arguments.
However, we see that with $L_I$ general, this will not be the case: 
the non-smooth locus of $X$ 
is strictly smaller-dimensional, so again using generality of $L_I$, we
see that every component of $p_I^{-1}(L_I) \cap X$ must have a dense
open subset inside the smooth locus of $X$. But $p_I^{-1}(L_I)$ is also
smooth, so we conclude that in this case, the intersection multiplicities
of every component of $p_I^{-1}(L_I) \cap X$ are simply determined by the
lengths of the intersected scheme, which is to say that
$[p_I^{-1}(L_I) \cap X] = [p_I^{-1}(L_I)] \cdot [X]$.
% irrelevant, but out of curiousity: can have imbedded comps in this case
% or no?

We also note that the proof of Proposition \ref{prop:hypersurface} shows
that if $X$ satisfies \eqref{eq:hypersurface-ineq}, and if
$p_I^{-1}(L_I) \cap X \neq \emptyset$ (still assuming $L_I$ general),
then we will have
that $p_I^{-1}(L_I) \cap X$ has generically finite fibers under 
$p_{I^c}$, so that 
$\dim p_{I^c} (p_I^{-1}(L_I) \cap X)=\dim p_I^{-1}(L_I) \cap X$.
Indeed, if $p_I^{-1}(L_I)$ is nonempty for a general $L_I$, this means
that $V_X$ maps dominantly to $\prod_{i \in I} {\rm Gr}(\alpha_i,n_i)$
under $p_{I} \circ \pi$, where $\pi$ is projection to the product of
Grassmannians. On the other hand, the proof of Proposition
\ref{prop:hypersurface} implies that there is a dense open subset $U \subseteq
V_X$ on which projection to $\prod_i {\rm Gr}(\alpha_i,n_i)$ is finite.
Thus, a general $L_I$ is in the image of $U$, meaning that there
exist $L_i$ for $i \in I^c$ such that $X \cap (L_1 \times \cdots \times L_k)$
is (nonempty and) finite. In particular, we must have that 
$p_I^{-1}(L_I) \cap X$ has finite fiber over $(L_i)_{i \in I^c}$, as
claimed.
% actually see U must be dense in general fiber, so every component has
% same dimension under p_{I^c}
\end{rem}

\begin{proof}[Proof of Proposition \ref{prop:form-evaluate}] 
Both sides being multiplicative, the desired identity
reduces immediately to the case that $\Xi=[X]$, with $X$ a subvariety.
Let $\widetilde{Y}=p_I^{-1}(L_I) \cap X$, let 
$\Psi=p_{I^c *}[\widetilde{Y}]$, and let 
$\bar{L}_I \in \prod_{i \in I} {\rm Gr}(\alpha_i,n_i)$ be the point
determined by the $L_i$.
Then on the level of underlying sets, we see that $Z_{\Psi,\bbeta_{I^c}}$ 
is given by the fiber of $Z_{X,\bbeta}$ over $\bar{L}_I$, which in turn is 
the vanishing cycle of partial evaluation of $H_{X,\beta}$ at $\bar{L}_I$.
Thus, we need 
to see that the associated multiplicities behave as expected. Let $V_X$ 
be the incidence correspondence in 
$\left(\prod_i \PP^{n_i}\right) \times \left(\prod_i {\rm Gr}(\alpha_i,n_i)\right)$,
and $V_{\Psi}$ be the incidence correspondence cycle as
in Remark \ref{rem:extension-compatible}, 
so that we have $Z_{[X],\bbeta}=\pi_* V_X$ and 
$Z_{\Psi,\bbeta_{I^c}}=\pi_* V_{\Psi}$ (although note that the
two maps $\pi$ are onto different products of Grassmannians).

Write
$\widetilde{L}_I \subseteq \prod_i {\rm Gr}(\alpha_i,n_i)$ 
for the fiber $p_I^{-1}(\bar{L}_I)$, and let 
$V_{\widetilde{Y}} \subseteq \left(\prod_i \PP^{n_i}\right) \times
\left(\prod_{i \in I^c} {\rm Gr}(\alpha_i,n_i)\right)$ 
be the scheme-theoretic incidence correspondence, which we will consider
as lying in 
$\left(\prod_i \PP^{n_i}\right) \times \left(\prod_i {\rm Gr}(\alpha_i,n_i)\right)$ 
using the point $\bar{L}_I$.  One then checks easily that 
$V_{\widetilde{Y}}=V_X \cap \pi^{-1} (\widetilde{L}_I)$,
for instance by comparing the functors of points.
% functors can be expressed in terms of tuples of line subbundles $\sL_i$
% of $\sO^{\oplus n_i+1}$ together with tuples of rk-$(\alpha_i+1)$ subbundles
% $\sV_i$ of same, with incidence condition that 
% $\sL_i \to \sO^{\oplus n_i+1}/\sV_i$ vanishes identically
% for first, assume $\sL_i$ induce point of $X$ and of $p_I^{-1}(L_I)$.
% impose incidence on I^c. extend out the $\sV_i$ for $i \in I$ using $L_i$.
% for second, assume $\sL_i$ induce point of $X$, and $\sV_i$ come from $L_i$
% for $i \in I$. impose incidence on all i. 
% for first, follows automatically after extending that have incidence for all
% i. for second, get point of $p_I^{-1}(L_I)$ by incidence and choice of 
% $\sV_i$ for i \in I. 
Next, we note that because of the generality of $L_I$, we have
$[V_X \cap \pi^{-1} (\widetilde{L}_I)]=[V_X] \cdot [\pi^{-1} (\widetilde{L}_I)]$.
Indeed, $V_X \cap \pi^{-1} \widetilde{L}_I$ is simply the fiber of $V_X$
over a general point of $\prod_{i \in I} {\rm Gr}(\alpha_i,n_i)$, and
the non-Cohen-Macaulay locus of $V_X$ is a proper algebraic subset, hence
of strictly smaller dimension. 
% proper subset quite general, but in this case can use generic smoothness
% since $V_X$ is a variety.
Semicontinuity of fiber dimension then implies that no component of a
general fiber is entirely contained in the non-Cohen-Macaulay locus of $V_X$,
and we obtain the desired identity as in Remark 
\ref{rem:intersect-agree}. The same argument shows that 
$[Z_{X,\bbeta} \cap \widetilde{L}_I]=[Z_{X,\bbeta}] \cdot [\widetilde{L}_I]$.

We next claim that 
$V_{\Psi}= (p_{I^c} \times p_{I^c})_*[V_{\widetilde{Y}}]$. 
This is clear again on the level of underlying sets, 
% direct from defns if do completely set-theoretically, replacing pushforward
% with image. can presumably show no comps have dim drops under $p_{I^c *}$,
% but OK even without this, since incidence correspondences are smooth over
% Y, so components and any dim drops line up on both sides.
so we just need to verify that the multiplicities
agree. By construction, $V_{\widetilde{Y}}$ is smooth over $\widetilde{Y}$,
so inherits the same multiplicities, and the two pushforwards under
$p_{I^c}$ visibly have the same fibers, so the claim follows.
% specifically, smoothness gives that comps and multiplicities of tilde{Y}
% agree with those of V_{tilde{Y}}. need to check degrees same for cycle
% pushforwards; for this enough to consider case tilde{Y} a variety. obvious
% fibers have same numbers of points, and scheme description above shows
% fibers actually scheme-theoretically isomorphic, so fine.

We then have
\begin{multline*} Z_{\Psi,\bbeta_{I^c}}=\pi_* V_{\Psi}
=\pi_*(p_{I^c} \times p_{I^c})_* [V_{\widetilde{Y}}]
= p_{I^c *} \pi_* [V_X \cap \pi^{-1}(\widetilde{L}_I)]
\\ = p_{I^c *} \pi_* ([V_X] \cdot [\pi^{-1}(\widetilde{L}_I)])
= p_{I^c *} ((\pi_* [V_X]) \cdot [\widetilde{L}_I]) 
\\ = p_{I^c *} ([Z_{[X],\bbeta}] \cdot [\widetilde{L}_I])
= p_{I^c *} [Z_{[X],\bbeta} \cap \widetilde{L}_I],\end{multline*}
% first and last equalities just defns, 2nd and 3rd verbatim from above,
% using associativity of pushforward. 5th should be projection formula,
% need works on levels of cycles (w/o ratl equivalence)
% 4th is important because naive projection can fail in non-CM case
where the fifth equality is the projection formula on the level of
cycles; see Proposition 8.1.1(c) of \cite{fu1}.
% here Y=Z, f=\pi, X'=|V_X|, Y'=|\pi(V_X)|, Z'=|\widetilde{L}_I|
Note that in the final expression, we are applying 
$p_{I^c *}$ to a cycle already supported in a fiber of $p_I$, so this
is just a formality, and we obtain the desired expression.
\end{proof}

We now conclude the desired assertion on multidegrees of
Cayley-Chow forms. We can extend multidegree linearly to cycles, so we
state the result in that context.

\begin{cor}\label{cor:multidegree}  
Under the hypotheses of Proposition \ref{prop:form-evaluate},
suppose that $\Xi$ has multidegree
$\sum_{\bgamma} 
a_{\bgamma} t_1^{\gamma_1} \cdots t_k^{\gamma_k}$. 
Then given $\bbeta$,
the Cayley-Chow form $H_{\Xi,\bbeta}$ has multidegree 
$$(a_{\alpha_1+1,\alpha_2,\dots,\alpha_k}, \dots, 
a_{\alpha_1, \dots, \alpha_{k-1},\alpha_k+1}).$$
\end{cor}

\begin{proof} By linearity, it suffices to treat the case that 
$\Xi=[X]$ for a subvariety $X$.
For each $j \in \{1,\dots,k\}$, 
we wish to show that the degree of $H_{X,\bbeta}$ in the 
$j$th set of variables is equal to 
\begin{equation}\label{eq:multideg}
a_{\alpha_1,\dots,\alpha_{j-1},\alpha_j+1,\alpha_{j+1},\dots,\alpha_k}.
\end{equation}
Setting $I=\{1,\dots,k\}\smallsetminus\{j\}$ and applying Proposition
\ref{prop:form-evaluate}, it suffices to show that the classical
Cayley-Chow form of
$p_{j*} (p_I^{-1}(L_I)\cap X))$ has degree given by \eqref{eq:multideg}
(note that it necessarily has the expected dimension; see Remark
\ref{rem:intersect-agree}). 
By the classical theory (Proposition 2.1 and 2.2 of \cite{g-k-z1}),
we thus want to show that $p_{j*} (p_I^{-1}(L_I)\cap X)$ has degree 
in $\PP^{n_j}$ given by \eqref{eq:multideg}. But this follows from the 
definitions and the projection formula in intersection theory (see for
instance p.\ 426 of \cite{ha1}).
% by defn, eq:multideg gives number of int pt w/p_j^{-1}(L_j') for genl
% L_j' of dim \alpha_j+1. proj formula says same as 
% number of pts of L_j' intersect $p_{j*} (p_I^{-1}(L_I)\cap X)$
\end{proof}

Finally, using the work of Castillo, Li and Zhang \cite{c-l-z1}, we can also 
translate the 
inequalities \eqref{eq:hypersurface-ineq} and \eqref{eq:good-ineqs} into 
multidegree-based criteria as follows. 

\begin{cor}\label{cor:multideg-ineqs}
If $X$ has multidegree $\sum a_{\bgamma} t_1^{\gamma_1} \cdots t_k^{\gamma_k}$, then the $Z_{X,\bbeta}$ associated to 
$(\beta_1,\dots,\beta_k)=(n_1-\alpha_1,\dots,n_k-\alpha_k)$ is a hypersurface 
if and only if 
\begin{equation}\label{eq:form-multideg}
a_{\alpha_1+1,\dots,\alpha_k} t_1 + \dots + a_{\alpha_1, \dots, \alpha_k+1} t_k
\end{equation}
is not identically zero. Moreover, $Z_{X,\bbeta}$ determines $X$ if and only if 
every term of \eqref{eq:form-multideg} is nonzero.
\end{cor}

\begin{proof} If some $n$-tuple $(\alpha_1,\dots,\alpha_j+1,\dots,\alpha_k)$ 
is in the support of the multi-degree, then Theorem~\ref{thm:support} 
implies that 
$|\bbeta_I| = \sum_{i \in I} (n_i-\alpha_i) \leq \dim p_I(X)+1$ 
for all $I$, so according to Proposition~\ref{prop:hypersurface} we have that
$Z_{X,\bbeta}$ is a hypersurface. 
Conversely, if 
$|\bbeta_I| \leq \dim p_I(X)+1$ for all $I$,
then we claim that there exists $j$ such that for any $I$ with
$$|\bbeta_I| = \dim p_I(X)+1,$$
we necessarily have $j \in I$.
Indeed, this follows immediately from Proposition \ref{prop:minimal-subset},
by choosing any $j \in J$. We then have
that $(\alpha_1,\dots,\alpha_j+1,\dots,\alpha_k)$ 
is in the support of the multidegree of $X$.

Next, if 
$|\bbeta_I| \leq \dim p_I(X)$ for all 
$I \subsetneq \{1,\dots,k\}$, it is clear that for all $j$ and $I$, we will 
have $\sum_{i \in I} (n_i - \gamma_i) \leq \dim p_{I}(X)$, where as before
$(\gamma_1,\dots, \gamma_k) = (\alpha_1,\dots,\alpha_j+1,\dots,\alpha_k)$.
Thus, for each $j$ we have $\bgamma$ in the support of the 
multidegree of $X$.
Conversely, if for some $I \subseteq \{1,\dots,k\}$ we have
$|\bbeta_I| > \dim p_I(X)$, then for any $j \not \in I$, with
$\bgamma$ as above we will have
$\sum_{i \in I} (n_i - \gamma_i) =|\beta_I| > \dim p_{I}(X)$,
so $\bgamma$ is not in the support of the multidegree of $X$. 
\end{proof}

\begin{rem}\label{rem:segre}
Given our basic setup, it is of course possible to re-embed $X$ into a
high-dimensional projective space via the Segre embedding, and then
apply the classical Cayley-Chow construction. This works canonically and
unconditionally to characterize $X$, but it doesn't reflect the geometry
of the embedding of $X$ into the original product of projective spaces, 
and it will typically require a great deal more data. For instance, if
a $3$-fold is embedded in $\PP^2 \times \PP^2$, then our construction
will give a bihomogeneous form in two sets of three variables. On the
other hand, the Segre embedding gives a $3$-fold in $\PP^8$, so the
relevent Grassmannian will be ${\rm Gr}(4,8)$, whose Pl\"ucker
embedding lands in $\PP^{125}$. Thus, the classical Cayley-Chow form is
in $126$ variables in this case! 
% plucker is in 9 choose 5 minus 1 dims.
\margh{work out degrees? say anything deeper?}
\end{rem}

We conclude with further examples. The first shows that in positive 
characteristic, our multigraded Cayley-Chow form may indeed come with
multiplicity strictly greater than $1$.

\begin{ex}\label{ex:frobenius} Let $K$ have characteristic $p$, and let
$X \subseteq \PP^2 \times \PP^2$ be the graph of the Frobenius morphism
$\varphi$.
Then $X$ has multidegree $p^2 t_1^2+pt_1t_2+ t_2^2$. If $\bbeta=(2,1)$, then
$$Z_{X,\bbeta}
=\{(P,L):\varphi(P) \in L\} \subseteq \PP^2 \times (\PP^2)^*.$$
If $L=Z(G)$ for a linear form $G$ and $P=(u_0,u_1,u_2)$, then 
$\varphi(P) \in L$ if and only if $G(u_0^p,u_1^p,u_2^p)=0$, so we see that
$Z_{X,\bbeta}$ is cut out by a $(p,1)$-form, as it should be. 
% Proposition \ref{prop:form-evaluate}: if plug in particular $u_0,u_1,u_2$,
% get linear form in coefs of $G$. Same as intersecting with preimage of
% point, taking image under $p_2$, and looking at lines containing it.
% or in other order, if fix $G$, get degree-$p$ form in coords of $P$.
% Intersecting with pullback of $G$ and taking image under $p_1$ gives
% p-thickened line, same thing. 

On the other hand, if $\bbeta=(1,2)$, then 
$$Z_{X,\bbeta}
=\{(L,P):P \in \varphi(L)\} \subseteq (\PP^2)^* \times \PP^2.$$
If $L=Z(G)$ and $P=(u_0,u_1,u_2)$ as above, then we observe that
$\varphi(L)$ is cut out by $\widehat{G}$, the linear form obtained from
$G$ by raising the coefficients to the $p$th power.
Then $P \in \varphi(L)$ if and only if $\widehat{G}(u_0,u_1,u_2)=0$, so
in this case $Z_{X,\bbeta}$ is still cut out by a $(p,1)$-form.
Thus, in order to get the right degree, we have to take the cycle 
$pZ_{X,\bbeta}$ in place of $Z_{X,\bbeta}$. We see this 
geometrically by observing that if we 
intersect $X$ with $p_2^{-1}(Y)$ for any point $Y$, we get a length-$p^2$
subscheme which can be identified under the first projection with
the fiber of $\vp$ over $Y$. Intersecting with a line in $\PP^2$ will
then reduce the length to $p$, but cannot reduce it to $1$. Thus,
the forgetful map from the incidence correspondence has degree $p$ in
this case.
% Proposition \ref{prop:form-evaluate}: if intersect with particular 
% $u_0,u_1,u_2$, get length $p^2$ thickening of $p$th roots of the $u_i$. 
% Then take $p^2$ power, get $\widehat{G}(u_0,u_1,u_2)^p$. Agrees with
% patched version of statement.
% if intersect with particular $L$ given by $G$, get variety whose image
% under $p_2$ is zero set of $\widehat{G}$. seems to just give
% $\widehat{G}(u_0,u_1,u_2)=0$? deeply confusing... ah, point is we
% shouldn't take image, rather take pushforward cycle. this should fix it,
% since induced map should be Frobenius on $\PP^1$, gives multiple of $p$.
\end{ex}

The following examples come from computer vision.

\begin{ex}\label{ex:multiview}[Multifocal tensors] We expand on the
discussion of multiview varieties from the introduction. Given $k \geq 2$
linear projections $\PP^3 \dashrightarrow \PP^2$ (viewed as positioned
pinhole cameras, with the centers of projection being camera centers), we 
let $X \subseteq (\PP^2)^k$ be the closure of the 
image of the induced rational map. This is called the ``multiview
variety'' associated to the camera configuration, and determines the
configuration (up to linear change of coordinates on $\PP^3$). If we 
assume that the
camera centers are all distinct, then $\dim p_I(X) = 3$ whenever $|I| \geq 2$, 
so Proposition~\ref{prop:hypersurface} guarantees that $Z_{X,\bbeta}$ 
is a hypersurface for all $(\beta_1,\dots, \beta_k)$ such that 
$|\bbeta| = 4$. 

On the other hand, $Z_{X,\bbeta}$ 
determines $X$ if and only if $|\bbeta_I| \leq 3$ for all 
$I \subsetneq \{1,\dots, k\}$, that is, if and only if $\beta_i \ne 0$ 
for all $i$. This clearly requires $k \leq 4$, and we see that if 
$k=2$ the vector $\bbeta$ must 
be $(2,2)$; if $k=3$ then it is a permutation of $(2,1,1)$; if $k=4$ then 
it is $(1,1,1,1)$. Moreover, the multidegree of $X$ is computed in
Corollary 3.5 of Aholt-Sturmfels-Thomas \cite{a-s-t1} to be
$$t_1^2\cdots t_k^2 \left(\sum_{1 \leq i_1 < i_2<i_3 \leq k} 
\frac{1}{t_{i_1}t_{i_2}t_{i_3}} + \sum_{1 \leq i_1,i_2 \leq k} 
\frac{1}{t_{i_1}^2 t_{i_2}}\right),$$
so according to Corollary~\ref{cor:multidegree}, we find that in the
allowed cases,
$H_{X,\bbeta}$ is multilinear, so it is associated with a tensor.

We thus recover the multifocal tensor construction when $k \leq 4$;
these are known as the ``fundamental matrix,'' the ``trifocal tensor,''
and the ``quadrifocal tensor,'' respectively.
On the other hand, for $k\geq 5$, we see that the constructed form never 
suffices to recover $X$. 
\margh{ also discuss need for nondegeneracy of camera centers?}
\end{ex}

\begin{ex}\label{ex:trifocal} Consider the $k=3$ (i.e., trifocal) case of 
the previous example, and let $P_i$ be the centers of projection.
It is helpful to observe that we
can think of $X$ as consisting of triples $(\ell_1,\ell_2,\ell_3)$ where
each $\ell_i$ is a line through $P_i$ in $\PP^3$, and 
$\ell_1 \cap \ell_2\cap\ell_3 \neq \emptyset$ (see Proposition 2.1 of
\cite{li5}). To avoid having to discuss
too many cases, we will assume that the $P_i$ are not collinear (and
in particular are distinct). As discussed in the previous example, we will 
have to have
up to permutation that $\bbeta=(2,1,1)$, so that in $\PP^2$,
we will have $L_1$ a point, and $L_2$ and $L_3$ lines. In the ambient
$\PP^3$, they will correspond to lines and planes containing $P_i$,
respectively. We will analyze the 
set $S_Z$ from the proof of Proposition \ref{prop:degenerate}. In fact,
the second author, Hebert and Ponce give a description of $S_Z$ in
Proposition 9 of \cite{t-h-p1}, observing that it does indeed contain
extra components beyond $X$ itself.

To describe $S_Z$, suppose we have fixed $(\ell_1,\ell_2,\ell_3)$,
not necessarily in $X$, so that we want to know under what conditions every 
$L_1,L_2,L_3$ containing $\ell_1,\ell_2,\ell_3$ must meet $X$, or 
equivalently, under what conditions every $L_2,L_3$ contain some choice
of $\ell_2',\ell_3'$ such that 
$\ell_1 \cap \ell_2' \cap \ell_3' \neq \emptyset$ (note that $L_1=\ell_1$
necessarily). Obviously, this is
the case if $\ell_1 \cap \ell_2 \cap \ell_3 \neq \emptyset$ already,
so that $(\ell_1,\ell_2,\ell_3) \in X$. However, there are two other
cases in which this occurs: if $\ell_1=\ell_2=\overline{P_1 P_2}$,
or $\ell_1=\ell_3=\overline{P_1 P_3}$.
Indeed, in the former case, we have that any plane $L_3$ must meet the line
$\ell_1=\ell_2$, yielding a choice of $\ell_3'$, and similarly for the 
latter case. On the other hand, one can check directly that in any other
situation, we can always find $L_2$, $L_3$ such that no $\ell_2',\ell_3'$
will have nonempty simultaneous intersection with $\ell_1$. Indeed, we
will be always be able to choose $L_2$ and $L_3$ so that $L_1 \cap L_2$
and $L_1 \cap L_3$ are distinct points, and $L_1 \cap L_2 \cap L_3$ is
therefore empty.
% if $L_1 \cap \ell_2$ is empty, can have $L_1 \cap L_2$ any point of $L_1$. 
% is $L_1 \cap \ell_2$ nonempty (nec not equal), general $L_2$ meets $L_1$
% at exactly that point. only bad case is both nonempty, and equal, which
% is case of being in $X$.
We thus find that $S_Z$ consists of $X$ together with two additional
$2$-dimensional components.

To compare this to the containment in \eqref{eq:degenerate}, we describe 
the geometry of the projections of $X$; by symmetry, it suffices to look 
at $p_1$ and $p_{\{1,2\}}$. 
We have that $p_1$ is surjective, and if we fix $\ell_1$, the fiber of $X$
over $\ell_1$ consists of pairs of lines $\ell_2,\ell_3$ which intersect
$\ell_1$ in a common point. If 
$\ell_1 \neq \overline{P_1 P_2},\overline{P_1 P_3}$, then each of 
$\ell_2$ and $\ell_3$ can meet $\ell_1$ in only a 
single point, so the choice of $\ell_2$ is determined by the choice of
a point on $\ell_1$, and then $\ell_3$ is determined as well. Thus, on
this set the fibers are $1$-dimensional. However, if 
$\ell_1=\overline{P_1 P_2}$ (so that $P_3 \not \in \ell_1$),
then every choice of $\ell_2$ meets $\ell_1$, and 
as long as $\ell_2 \neq \ell_1$, then $\ell_1 \cap \ell_2 = \{P_2\}$, and
$\ell_3=\overline{P_2 P_3}$ is uniquely determined. On the other hand, we 
could also have
$\ell_2=\ell_1$, in which case we have a $1$-dimensional set of choices of
$\ell_3$. Thus, in this case the fiber has two components, one of dimension
$2$, and one of dimension $1$. The same holds if $\ell_1=\overline{P_1 P_3}$. 
We conclude that the general fiber is $1$-dimensional, but there
are two fibers which are (non-purely) $2$-dimensional, corresponding to 
$\overline{P_1 P_2}$ and $\overline{P_1 P_3}$, respectively. The analogous 
description holds for $p_2$ and $p_3$.

Next, the image of $p_{\{1,2\}}$ is precisely the set of pairs $(\ell_1,\ell_2)$
which have nonempty intersection, which forms a $3$-dimensional set. 
Provided $\ell_1 \neq \ell_2$, we will have that $\ell_1 \cap \ell_2$ is
a single point. If this point is not $P_3$, then $\ell_3$ is uniquely
determined, and thus $p_{\{1,2\}}$ is injective over such pairs. If 
$\ell_1 =\ell_2=\overline{P_1 P_2}$, so that $P_3$ is not on $\ell_1$
or $\ell_2$,
then $\ell_3$ is determined by a choice of point of $\ell_1$, and we
have a $1$-dimensional fiber. Finally, if $\ell_1 \neq \ell_2$ but
both go through $P_3$ (so that they are necessarily $\overline{P_1 P_3}$
and $\overline{P_2 P_3}$ respectively), then any choice of $\ell_3$ is
valid, and we obtain a $2$-dimensional fiber. To summarize, the general
fiber is $0$-dimensional, but the fiber corresponding to 
$(\overline{P_1 P_2},\overline{P_1 P_2})$ is $1$-dimensional, and the fiber
corresponding to 
$(\overline{P_1 P_3},\overline{P_2 P_3})$ 
is $2$-dimensional. The analogous description holds for $p_{\{2,3\}}$ and
$p_{\{1,3\}}$.

We compare this to the proof of Proposition \ref{prop:degenerate} as 
follows: if $I=\{i\}$, 
the set $X_I$ is where the fibers of $p_i$ have dimension at least
$2$ if $i=1$, and at least $3$ if $i=2,3$. The latter two cases give
the empty set, but the former consists of the two points where $\ell_1$
is either $\overline{P_1 P_2}$ or $\overline{P_1 P_3}$. Note that
the two extra components of $S_Z$ we have identified are contained in
$p_1^{-1}(X_{\{1\}})$, but strictly. On the other hand, if $I=\{i,j\}$,
the set $X_I$ is where the fibers of $p_{\{i,j\}}$ have dimension at
least $2$ if $I=\{2,3\}$, and at least $1$ otherwise. In the first case,
we get that $X_I$ is the single point where $\ell_2=\overline{P_1 P_2}$
and $\ell_3=\overline{P_1 P_3}$. In this case, we already have that
$p_{\{2,3\}}^{-1}(\ell_2,\ell_3)$ is contained in $X$, so we do not get any
new component of $S_Z$. For $I=\{1,2\}$, we have the same behavior over
$(\overline{P_1 P_3},\overline{P_2 P_3})$, but $X_I$ also includes the
point $(\overline{P_1 P_2},\overline{P_1 P_2})$, and one additional 
component of $S_Z$ is equal to the fiber of $p_{\{1,2\}}$ over this point. 
Considering $I=\{1,3\}$ gives the other additional component of $S_Z$.
Thus, in this case we have strict containment in \eqref{eq:degenerate},
although we obtain equality if we restrict $I$ to two-element sets.
\margh{odd: collinear affects geometry, but doesn't screw up tensor? even
when two points coincide, may have one okay choice of tensor?}
\end{ex}

\bibliographystyle{amsalpha}
\bibliography{gen}

\end{document}